\documentclass[3p,times]{article}
\usepackage{bm,ulem}

\newcommand{\R}{\mathbf{R}}
\newcommand{\B}{\mathbf{B}}
\renewcommand{\P}{\mathbf{P}}

\newcommand{\x}{\mathbf{x}}
\newcommand{\y}{\mathbf{y}}

\newcommand{\M}{\mathbf{M}}

\newcommand{\J}{\mathcal{J}}
\newcommand{\A}{\mathbf{a}}

\newcommand{\nvar}{\textsc{n}_\textnormal{var}}

\usepackage{bbm}
\usepackage{subfigure}
\usepackage{amsmath}
\usepackage{amssymb}
\usepackage{graphicx}
\usepackage{capt-of}
\usepackage{color}
\usepackage{url}
\usepackage{algorithm}
\usepackage{algpseudocode,caption}
\usepackage{mdwlist}
\usepackage{array}
\usepackage{tikz}
\usepackage{dsfont}
\usepackage{multirow}
\usepackage[toc,page]{appendix}

\algblock{ParFor}{EndParFor}
\algnewcommand\algorithmicparfor{\textbf{For all}}
\algnewcommand\algorithmicpardo{\textbf{do in parallel}}
\algnewcommand\algorithmicendparfor{\textbf{end\ For all}}
\algrenewtext{ParFor}[1]{\algorithmicparfor\ #1\ \algorithmicpardo}
\algrenewtext{EndParFor}{\algorithmicendparfor}

\title{A Time-parallel Approach to Strong-constraint Four-dimensional Variational Data Assimilation}
\author{
	Vishwas Rao and Adrian Sandu
}
\date{\today}
\begin{document}
 \thispagestyle{empty}
\setcounter{page}{0}

\begin{Huge}
\begin{center}
Computational Science Laboratory Technical Report CSL-TR-18-2015\\
\today
\end{center}
\end{Huge}
\vfil
\begin{huge}
\begin{center}
Vishwas Rao and Adrian Sandu
\end{center}
\end{huge}

\vfil
\begin{huge}
\begin{it}
\begin{center}
``A Time-parallel Approach to Strong-constraint Four-dimensional Variational Data Assimilation''
\end{center}
\end{it}
\end{huge}
\vfil

\begin{large}
\begin{center}
Computational Science Laboratory \\
Computer Science Department \\
Virginia Polytechnic Institute and State University \\
Blacksburg, VA 24060 \\
Phone: (540)-231-2193 \\
Fax: (540)-231-6075 \\  
Email: \url{sandu@cs.vt.edu} \\
Web: \url{http://csl.cs.vt.edu}
\end{center}
\end{large}

\vspace*{1cm}

\begin{tabular}{ccc}
\includegraphics[width=2.5in]{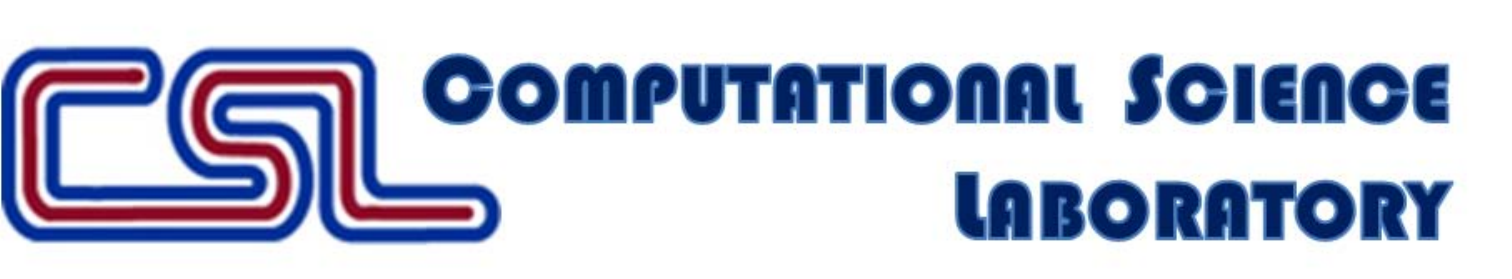}
&\hspace{2.5in}&
\includegraphics[width=2.5in]{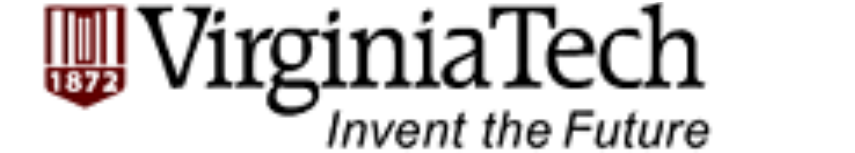} \\
{\bf\em Innovative Computational Solutions} &&\\
\end{tabular}

\newpage

 \maketitle
\begin{abstract}
A parallel-in-time algorithm based on an augmented Lagrangian approach is proposed to solve four-dimensional variational (4D-Var) data assimilation problems. The assimilation window is divided into multiple sub-intervals that allows to parallelize cost function and gradient computations. Solution continuity equations across interval boundaries are added as constraints. The augmented Lagrangian approach leads to a different formulation of the variational data assimilation problem than weakly constrained 4D-Var.  A combination of serial and parallel 4D-Vars to increase performance is also explored. The methodology is illustrated on data assimilation problems with Lorenz-96 and the shallow water models. 
\end{abstract}
\newpage
\tableofcontents

\newpage\setcounter{page}{1}
%

\section{Introduction} \label{sec:Intro}
Predicting the behavior of complex dynamical systems, such as the atmosphere, requires  using information from observations to decrease the uncertainty in the forecast. Data assimilation combines information from a numerical model, prior knowledge, and observations (all with associated errors) in order to obtain an improved estimate of the true state of the system. Data assimilation is an important application of data-driven application systems (DDDAS \cite{dddas2}, or InfoSymbiotic systems) where measurements of the physical system are used to constrain simulation results.

Two approaches to data assimilation have gained widespread popularity: variational  and ensemble-based methods. The ensemble-based methods are rooted in statistical theory, whereas the variational approach is derived from optimal control theory. The variational approach formulates data assimilation as a nonlinear optimization problem constrained by a numerical model. The initial conditions (as well as boundary conditions, forcing, or model parameters) are adjusted to minimize the discrepancy between the model trajectory and a set of time-distributed observations. In real-time operational settings the data assimilation process is performed in cycles: observations within an assimilation window are used to obtain an optimal trajectory, which provides the initial condition for the next time window, and the process is repeated in the subsequent cycles. The variational methodology is widely adopted by most national and international numerical weather forecast centers to provide the initial state for their forecast models. 
       
Performing nonlinear optimization in the 4D-Var framework is an inherently sequential process. Computer architectures  progressively incorporate more parallelism, while maintaining a constant processor speed. As our understanding of the physics improves the computer models become increasingly more complex. Advanced and scalable parallel algorithms to solve 4D-Var need to be developed to continue to perform data assimilation in real time. This challenge has been partially addressed by exploring parallelism in spatial dimension. Tr\'emolet and Le Dimet  \cite{Tremolet1996657} have shown how variational data assimilation can be used to couple models and to perform parallelization in space for the assimilation process. Rantakokko \cite{Rantakokko19972017} considers different data distribution strategies to perform parallel variational data assimilation in the spatial dimension. A scalable approach for three dimensional variational (3D-Var) data assimilation is presented in \cite{Damore:2014} and parallelism is achieved by dividing the global problem into multiple local 3D-Var sub-problems. Multiple copies of modified 3D-Var problem, which ensures feasibility at the boundaries of the sub-domains, are solved across processors and the global 3D-Var minimum is obtained by collecting the local minima.
       
An important challenge associated with computing solutions to 4D-Var problem is parallelization in the temporal dimension. Fisher  \cite{fisher2013Parallel4DVar} attempts to address this challenge  by the saddle point formulation that solves directly the optimality equations. Aguiar et al. \cite{ICML2011Martins_150} apply the augmented Lagrangian approach to constrained inference problems in the context of graphical models. The approach proposed herein uses  the augmented Lagrangian framework in the context of 4D-Var data assimilation. The most computationally expensive components of 4D-Var, namely the cost function and gradient evaluations, are performed in a time-parallel manner.
       
The remainder of this paper is organized as follows: Section \ref{sec:DA} introduces data assimilation and the traditional 4D-Var approach. Section \ref{sec:Reformulation} formulates the 4D-Var problem in an augmented Lagrangian framework to expose the time parallelism in the cost function and gradient evaluations. Section \ref{sec:Algorithm} gives a detailed description of the parallel assimilation algorithm. Section \ref{sec:exp} shows the numerical results with the small, chaotic Lorenz-96 model, and a relatively large shallow water on the sphere model. Concluding remarks and future research directions are discussed in Section \ref{sec:conc}.

\section{Four-dimensional variational data assimilation} \label{sec:DA}

Data assimilation (DA) is the fusion of information from priors, imperfect model predictions, and noisy data, to obtain a consistent description of the true state $\x^{\rm true}$ of a physical system \cite{daley1993, kalnay2003, sandu2011chemical, sandu2005adjoint}. The best estimate that optimally fuses all these sources of information is called the analysis $\x^{\rm a}$. 

The prior information encapsulates our current knowledge of the system. Usually the prior information is contained in a background estimate of the state $\x^{\rm b}$ and the corresponding background error covariance matrix $\mathbf{B}$. 

The model captures our knowledge about the physical laws that govern the evolution of the system. The model evolves an initial state $\x_0 \in \mathbb{R}^n$ at the initial time $t_0$ to future states $\x_{\rm k} \in \mathbb{R}^n$ at future times $t^{\{ \ell \}}$. A general model equation is represented as follows:
\begin{equation}
\label{eqn:genmodel}
 \x_{\rm k} = \mathcal{M}_{t_0 \rightarrow t^{\{ \ell \}}} \left(\x_0\right)\,.
\end{equation}

Observations are noisy snapshots of reality available at discrete time instances. Specifically, measurements $\y_{\rm k} \in \mathbb{R}^m$ of the physical state $\x^{\rm true}\left(t^{\{ \ell \}}\right)$ are taken at times $t^{\{ \ell \}}$, ${\rm k}=1,\cdots, N$.
%
The model state is related to observations by the following relation:
\begin{eqnarray}
\label{eqn:genobs}
 \y_{\rm k} &=& \mathcal{H}\left(\x_{\rm k}\right)-\varepsilon_{\rm k}^{\rm obs},\quad {\rm k}=1,\cdots,N, \\
 \nonumber
 \varepsilon_{\rm k}^{\rm obs} &=& \varepsilon_{\rm k}^{\rm representativeness} + \varepsilon_{\rm k}^{\rm measurement} \,.
\end{eqnarray}
The observation operator $\mathcal{H}$ maps the model state space onto the observation space. The observation error term $\left(\varepsilon_{\rm k}^{\rm obs}\right)$ accounts for both measurement and representativeness errors. Measurement errors are due to imperfect sensors. The representativeness errors are due to the inaccuracies of the mathematical and numerical approximations inherent to the model.

Variational methods solve the data assimilation problem in an optimal control framework, where one finds the control variable which minimizes the mismatch between the model forecasts and the observations. Strong-constraint 4D-Var assumes that the model  \eqref{eqn:genmodel} is perfect \cite{sandu2011chemical, sandu2005adjoint}. The control parameters are the initial conditions $\x_0$, which uniquely determine the state of the system at all future times via the model equation \eqref{eqn:genmodel}. The background state is the prior best estimate of the initial conditions $\x_0^{\rm b}$, and has an associated initial background error covariance matrix $\mathbf{B}_0$. Observations $\y_{\rm k}$ at $t^{\{ \ell \}}$ have the corresponding observation error covariance matrices $\mathbf{R}_{\rm k}$, ${\rm k}=1,\cdots,N$. The 4D-Var problem provides the estimate $\x_0^{\rm a}$ of the true initial conditions as the solution of the following optimization problem
\begin{equation}\label{eqn:ip}
 \x_0^{\rm a}  =  \underset{\x_0} {\text{ arg\, min}}~~ \J\left(\x_0\right) \qquad
 \text{subject to}~ \text{\eqref{eqn:genmodel}},
\end{equation}
with the following cost function:
\begin{eqnarray}\label{eqn:fdvar-cf}
 \mathcal{J}\left(\x_0\right) &=& \frac{1}{2} \left(\x_0 - \x_0^{\rm b} \right)^{\rm T} \mathbf{B}_0^{-1} \left(\x_0 - \x_0^{\rm b} \right) \\ \nonumber
 &&+ \frac{1}{2} \displaystyle \sum_{\rm k=1}^N \left(\mathcal{H}_{\rm k}\left(\x_{\rm k}\right) - \y_{\rm k} \right)^{\rm T} \mathbf{R}_{\rm k}^{-1} \left(\mathcal{H}_{\rm k}\left(\x_{\rm k}\right) - \y_{\rm k} \right).
\end{eqnarray}
The first term of the sum \eqref{eqn:fdvar-cf} quantifies the departure of the solution $\x_0$ from the background state $\x_0^{\rm b}$ at the initial time $t_0$. The second term measures the mismatch between the forecast trajectory (model solutions $\x_{\rm k}$) and observations $\y_{\rm k}$ at all times $t^{\{ \ell \}}$ in the assimilation window. The weighting matrices $\mathbf{B}_0$ and $\mathbf{R}_{\rm k}$ need to be predefined, and their quality influences the accuracy of the resulting analysis.

Weak constraint 4D-Var \cite{sandu2011chemical} removes the perfect model assumption by allowing a model error 
$\eta_{{\rm k}+1} = \x_{\rm k +1} - \mathcal{M}_{\rm k, \rm k+1} \left(\x_{\rm k}\right)$. Under the assumption that the model errors are normally distributed, $\eta_{\rm k} \in \mathcal{N}(0, \mathbf{Q}_{\rm k})$, the weak constraint 4D-Var solution is the unconstrained minimizer of the cost function
\begin{eqnarray}\label{eqn:weakcf}
 \mathcal{J}^{\rm weak}\left(\x_0,\dots,\x_N\right) &=& \mathcal{J}\left(\x_0\right) + \\\nonumber
 &&\frac{1}{2} \displaystyle \sum_{k=0}^{N-1} \left(\x_{\rm k +1} - \mathcal{M}_{\rm k, \rm k+1} \left(\x_{\rm k}\right)\right)^{\rm T} \mathbf{Q}_{{\rm k}+1}^{-1} \, \left(\x_{\rm k +1} - \mathcal{M}_{\rm k, \rm k+1} \left(\x_{\rm k}\right)\right).
\end{eqnarray}
The control variables are the states of the system at all times in the assimilation window. 

In this paper we focus on the strong constraint formulation \eqref{eqn:ip}. The minimizer of \eqref{eqn:ip} is computed iteratively using gradient-based numerical optimization methods. First-order adjoint models provide the gradient of the cost function \cite{cacuci2005sensitivity}, while second-order adjoint models provide the Hessian-vector product (e.g., for Newton-type methods). The methodology for building and using various adjoint models for optimization, sensitivity analysis, and uncertainty quantification is discussed in \cite{cioaca2012second, sandu2005adjoint}. Various strategies to improve the the 4D-Var data assimilation system are described in \cite{cioaca2014optimization}. The procedure to estimate the impact of observation and model errors is developed in \cite{rao2014posterioriJournal, rao2014posteriori}. A framework to perform derivative free variational data assimilation using the trust-region framework is given in \cite{ruiz2015derivative}.

The iterative solution of \eqref{eqn:ip} is highly sequential: first, one iteration follows the other; next, the forward and adjoint models are run sequentially forward and backward in time, respectively.  In order to reveal additional parallelism the solution to 4D-Var problem is is approached using the augmented Lagrangian framework. In this framework, the assimilation window is divided into multiple sub-intervals, and the model constraints are explicitly imposed at the boundaries. This approach bears similarities with the Parareal approach that exploits time parallelism in the solution of ordinary differential equations \cite{Rao201476}. 

\section{4D-Var solution by the augmented Lagrangian approach} \label{sec:Reformulation}
The 4D-Var cost function \eqref{eqn:fdvar-cf} is minimized subject to the generic model constraints \eqref{eqn:genmodel}. The model equations can also be written as
\begin{equation}
\label{eqn:model}
 \x_{\rm k +1} = \mathcal{M}_{\rm k, \rm k+1} \left(\x_{\rm k}\right), \quad {\rm k} = 0,1,\cdots,N-1\,,
\end{equation}
where $\mathcal{M}_{\rm k, \rm k+1}$ represents the model solution operator that propagates the state $\x_{k}$ at $t_{k}$ to the state $\x_{k+1}$ at $t_{k+1}$. 
The minimizer of  \eqref{eqn:fdvar-cf} under the constraints \eqref{eqn:model} is the unconstrained minimizer of the Lagrangian 
\begin{eqnarray}
\label{eqn:Lagrangian}
 \mathbf{L}(\x_0,\dots,\x_N;\lambda_0,\dots,\lambda_N) &=& 
  \mathcal{J}(\x_0)-\displaystyle \sum_{\rm k=0}^{\rm N-1} \lambda_{{\rm k}+1}^{\rm T} \cdot \left(\x_{\rm k +1} - \mathcal{M}_{\rm k, \rm k+1} \left(\x_{\rm k}\right)\right)\,.
\end{eqnarray}
To expose time parallelism the assimilation window is divided into $N$ sub-intervals, namely, 
\begin{equation}
\label{eqn:interval-partition}
\lbrack t_0, t_N \rbrack = \lbrack t_0, t_1\rbrack \, \cup \,\dots \, \cup \, \lbrack t_{N-1}, t_N\rbrack.
\end{equation}
The forward model and adjoint model states at the interval boundaries are denoted by 
\begin{equation}
\label{eqn:solution-partition}
\x = [\x_0,\, \cdots, \,\x_N], \quad \bm{\lambda} = [\lambda_0,\, \cdots, \,\lambda_N],
\end{equation}
respectively. We denote by $\x^{\rm a}$ the optimal solution and by $\bm{\lambda}^{\rm a}$ the optimal value of the Lagrange multiplier in \eqref{eqn:Lagrangian}.

The augmented Lagrangian  \cite[Section 17.3]{wright1999numerical} associated with \eqref{eqn:fdvar-cf} and \eqref{eqn:model} reads:
\begin{eqnarray}
\label{eqn:AugLagCF}
 \mathcal{L}\left(\x,\bm\lambda,\mu\right) &=& \frac{1}{2} \left(\x_0 - \x_0^{\rm b} \right)^{\rm T} \mathbf{B}_0^{-1} \left(\x_0 - \x_0^{\rm b} \right) \\ \nonumber
 &&+ \frac{1}{2} \sum_{\rm k=1}^{\rm N}\left(\mathcal{H} \left(\x_{\rm k}\right) - \y_{\rm k}\right)^{\rm T} \R_{\rm k}^{-1} \left(\mathcal{H} \left(\x_{\rm k}\right) - \y_{\rm k}\right)\\ 
 \nonumber
 && -\sum_{\rm k=0}^{\rm N-1} \lambda_{\rm k+1}^{\rm T} \left(\x_{\rm k+1} -\mathcal{M}_{\rm k, k+1} \left(\x_{\rm k}\right)\right) \\\nonumber 
 && + \frac{\mu}{2} \sum_{\rm k=0}^{\rm N-1} \left(\x_{\rm k+1} - \mathcal{M}_{\rm k, k+1} \left(\x_{\rm k}\right)\right)^{\rm T}\, \mathbf{P}_{\rm k+1}^{-1}\,\left(\x_{\rm k+1} - \mathcal{M}_{\rm k, k+1} \left(\x_{\rm k}\right)\right)\,,
\end{eqnarray}
where $\P_{\rm k}$'s are error scaling matrices. This is the (regular) Lagrangian for the problem that minimizes the cost function
\begin{eqnarray}
\label{eqn:augmentedcf}
 \mathcal{J}^{\rm augmented}\left(\x_0,\dots,\x_N\right) &=& \mathcal{J}\left(\x_0\right) + \\ \nonumber
 &&\frac{\mu}{2} \displaystyle \sum_{k=0}^{N-1} \left(\x_{\rm k +1} - \mathcal{M}_{\rm k, \rm k+1} \left(\x_{\rm k}\right)\right)^{\rm T} \mathbf{P}_{{\rm k}+1}^{-1} \, \left(\x_{\rm k +1} - \mathcal{M}_{\rm k, \rm k+1} \left(\x_{\rm k}\right)\right)
\end{eqnarray}
subject to the model constraint \eqref{eqn:model}. The constrained minimization of \eqref{eqn:augmentedcf} is equivalent to\eqref{eqn:ip} since the additional term is zero along the constraints. Note that \eqref{eqn:augmentedcf} is a constrained minimization problem, unlike \eqref{eqn:weakcf}, which is unconstrained.

The original 4D-Var problem in \eqref{eqn:ip} is solved in the augmented Lagrangian framework by performing a sequence of unconstrained minimizations
\begin{equation}
\label{eqn:ipAug}
 \x^{\{\ell\}}  =  \arg\, \min_{\x}~~ \mathcal{L}\left(\x, \widetilde{\bm\lambda}^{\{\ell\}}, \mu^{\{\ell\}} \right), \quad \ell = 0,1,\dots\,. 
\end{equation}
If $\widetilde{\bm\lambda}^{\{\ell\}} \approx \bm\lambda^{\rm a}$ then $\widetilde{\x}^{\{\ell\}} \approx \x^{\rm a}$, and the solution error decreases with increasing $\mu$ \cite[Section 17.3]{wright1999numerical}. 

The optimization proceeds in cycles of inner and outer iterations. Inner iterations solve the optimization problem \eqref{eqn:ipAug} for particular values of $\mu^{\{\ell\}}$ and $\widetilde{\bm\lambda}^{\{\ell\}}$. After each solution of \eqref{eqn:ipAug} the outer iteration $\ell$ is completed by
updating the Lagrange multiplier approximation and the penalty parameter, as follows:
\begin{subequations}\label{eqn:update}
\begin{eqnarray}
\label{eqn:updateMu}
\mu^{\rm \{\ell+1\}} &=& \rho\, \mu^{\rm  \{\ell\}}\,, \\[3pt]
\label{eqn:updateLambda}
\lambda^{\rm \{\ell+1\}}_{\rm k} &=& \lambda^{\rm \{\ell\}}_{\rm k} - \mu^{\rm \{\ell\}} \left(\x^{\rm \{\ell\}}_{\rm k} - \mathcal{M}_{\rm k-1, k}(\x^{\rm \{\ell\}}_{\rm k-1}) \right)\,, \quad \rm{k}=0, \dots , N\,.
\end{eqnarray}
\end{subequations}
The penalty parameter is progressively increased by a constant $\rho > 1$ in order to impose the model constraints \eqref{eqn:model}. Different other strategies to update $\mu$ and $\bm\lambda$ can be used \cite{birgin2008improving}.

\subsection{Augmented Lagrangian optimization}

Figure \ref{fig:Iterations_Lorenz} illustrates the convergence process of the augmented Lagrangian 4D-Var. The assimilation window is divided into multiple sub-intervals \eqref{eqn:interval-partition}. The control vector $\x$ for the optimization process contains the state vector dat the beginning of each of the sub-intervals \eqref{eqn:interval-partition}. The initial value of the control vector contains the background states $\x^{\rm b}_k$ at the beginning of each of the sub-interval, and is therefore a continuous curve.  The outer iterations start with a small value of $\mu$ that only imposes the constraints \eqref{eqn:model} loosely. Consequently, the solutions during the first outer iterations show large discontinuities at the interval boundaries. Subsequently, $\mu$ is increased and as a result the constraints are satisfied accurately. This results in a smooth solution curve resembling closely the serial 4D-Var solution.

\section{The parallel algorithm}\label{sec:Algorithm}

Most of the computational time required by the 4D-Var  solution is claimed by the multiple cost function and gradient evaluations. In the augmented Lagrangian 4D-Var formulation \eqref{eqn:ipAug}--\eqref{eqn:update} the forward and adjoint models can be run in parallel over sub-intervals, as explained next.

\subsection{Parallel-in-time runs of the forward model}

The value of the augmented Lagrangian cost function  \eqref{eqn:AugLagCF} can be computed in parallel. Specifically, on each sub-interval $[t^{\{ \ell \}},t^{\{ \ell+1 \}}]$ in \eqref{eqn:interval-partition} a forward solution $\mathcal{M}_{\rm k, k+1} \left(\x_{\rm k}\right)$ is computed starting from the initial value $\x_{\rm k}$  \eqref{eqn:solution-partition}. The forward model runs on each sub-interval can be carried out concurrently. The computational steps are detailed in Algorithm \ref{alg:cf}.
\begin{algorithm}
\caption{Parallel\_CostFunction }\label{alg:cf}
\begin{algorithmic}[1]
\Procedure{Parallel\_CostFunction}{}
\State \textbf{Input:} $\lbrack \x_0, \dots, \x_{\rm N};\, \lambda_1, \dots, \lambda_{\rm N} \rbrack$
\State \textbf{Output:} $\mathcal{L}$
\State $\mathcal{L} \gets 0$
\ParFor{$0\leq {\rm k} \leq N-1$}
\State $\Delta\x_{\rm k+1} \gets \,\x_{\rm k+1} - \mathcal{M}_{\rm k,k+1}\left(\x_{\rm k}\right)$   \Comment{\parbox[t]{.3\linewidth}{\small{Solution mismatch at sub-interval boundaries}}}
\State $\Delta \y_{\rm k+1} \gets \mathcal{H} \left(\x_{\rm k+1}\right) - \y_{\rm k+1}$
\Comment{\parbox[t]{.3\linewidth}{\small{Solution mismatch with observations}}}
\EndParFor 
\For{$0\leq {\rm k} \leq N-1$} \Comment{\parbox[t]{.3\linewidth}{\small{Evaluate the cost function using \eqref{eqn:AugLagCF}}}}
\State $\mathcal{L} \gets \mathcal{L} + \displaystyle\frac{\mu}{2}\, \Delta \x_{\rm k+1}^{\rm T}\, \mathbf{P}^{-1}_{\rm k+1}\, \Delta \x_{\rm k+1}
 -  \lambda_{\rm k+1}^{\rm T}\, \Delta \x_{\rm k+1} + \displaystyle\frac{1}{2}\, \Delta \y_{\rm k+1}\,  \R_{\rm k+1}^{-1}\, \Delta \y_{\rm k+1}$
\EndFor
\State $\mathcal{L} \gets \mathcal{L}+\displaystyle\frac{1}{2} \left(\x_0 - \x_0^{\rm b} \right)^{\rm T} \mathbf{B}_0^{-1} \left(\x_0 - \x_0^{\rm b} \right)$
\EndProcedure
\end{algorithmic}
\end{algorithm}

\subsection{Parallel-in-time runs of the adjoint model}

The gradient of the augmented Lagrangian \eqref{eqn:AugLagCF} with respect to the forward model state is given by:
\begin{subequations}
\label{eqn:AugLagGrad}
\begin{eqnarray}
\label{eqn:AugLagGrad0}
\nabla_{\x_0} \mathcal{L}&=&\B_0^{-1} \left(\x_0 - \x_0^{\rm b}\right) + \M^{\rm T}_{0,1} \lambda_1 - {\mu}\, \M^{\rm T}_{0,1} \, \P_0^{-1} \left(\x_1 - \mathcal{M}_{0,1}\left(\x_0\right) \right)\,, \\ [5pt]
\label{eqn:AugLagGrad1} 
\nabla_{\x_{\rm k}}\mathcal{L} &=&\mathbf{H}_{\rm k}^{\rm T} \R_{\rm k}^{-1} \left(\mathcal{H}\left(\x_{\rm k}\right) -\y_{\rm k} \right) - \lambda^{\rm T}_{\rm k} + \mathbf{M}^{\rm T}_{\rm k,\rm k+1} \lambda_{\rm k+1} \\
&& + \mu\, \mathbf{P}^{-1}_{{\rm k}-1}\left(\x_{\rm k} - \mathcal{M}_{\rm k-1,k}\left(\x_{\rm k-1}\right)\right)  \nonumber  \\ 
 && - \mu\, \mathbf{M}^{\rm T}_{\rm k,\rm k+1}\, \mathbf{P}^{-1}_{\rm k} \left(\x_{\rm k+1} - \mathcal{M}_{\rm k, k+1} \left(\x_{\rm k}\right)\right)\,, \nonumber  
 \nonumber
  \quad \ {\rm k} = 1, \dots,  N-1, \\ [5pt]
\label{eqn:AugLagGrad2} 
\nabla_{\x_{\rm N}}\mathcal{L} &=& \mathbf{H}_{\rm N}^{\rm T} \R_{\rm N}^{-1}\left(\mathcal{H}\left(\x_{\rm N}\right) -\y_{\rm N} \right) - \lambda_{\rm N}^{\rm T}\\ \nonumber
&&+ \mu\, \P_{\rm N -1}^{-1} \left( \x_{\rm N} - \mathcal{M}_{\rm{N}-1, \rm N}\left(\x_{\rm N-1}\right) \right)\,,
\end{eqnarray}
\end{subequations}
where the tangent linear operators of the observation and model solution operators are
\[
\left.\mathbf{H}_{\rm k} = \frac{\partial \mathcal{H}_{\rm k}}{\partial \x_{\rm k}}\right\rvert_{\x_{\rm k}},
\qquad \left. \mathbf{M}_{\rm k,k+1} = \frac{\partial \mathcal{M}_{\rm k,k+1}}{\partial \x_{\rm k}}\right\rvert_{\x_{\rm k}},
\]
respectively, and their transposes are the corresponding adjoint operators. Using the notation
\begin{eqnarray*}
\Delta\x_{\rm k} &:=& \x_{\rm k} - \mathcal{M}_{\rm k-1, k} \left(\x_{\rm k-1}\right) \\[5pt]
\Delta\y_{\rm k} &:=& \mathcal{H}\left(\x_{\rm k}\right) -\y_{\rm k}
\end{eqnarray*}
the gradient \eqref{eqn:AugLagGrad} can be written as
\begin{eqnarray*}
\nabla_{\x_0} \mathcal{L}&=&\B_0^{-1} \left(\x_0 - \x_0^{\rm b}\right) + \M^{\rm T}_{0,1} \lambda_1 - {\mu}\, \M^{\rm T}_{0,1} \, \P_1^{-1} \Delta\x_1\,, \\ [5pt]
\nabla_{\x_{\rm k}}\mathcal{L} &=&\mathbf{H}_{\rm k}^{\rm T} \R_{\rm k}^{-1}\, \Delta\y_{\rm k} - \lambda^{\rm T}_{\rm k} + \mathbf{M}^{\rm T}_{\rm k,\rm k+1}\, \left( \lambda_{\rm k+1} - \mu \, \mathbf{P}^{-1}_{\rm k+1} \Delta\x_{\rm k+1}\right) \\
&& + \mu\, \mathbf{P}^{-1}_{{\rm k}}\, \Delta\x_{\rm k},  \nonumber  
 \nonumber
  \quad \ {\rm k} = 1, \dots,  N-1, \\ [5pt]
\nabla_{\x_{\rm N}}\mathcal{L} &=& \mathbf{H}_{\rm N}^{\rm T} \R_{\rm N}^{-1}\, \Delta\y_{\rm N} - \lambda_{\rm N}^{\rm T} + \mu\, \P_{\rm N -1}^{-1} \, \Delta\x_{\rm N}\,.
\end{eqnarray*}

The augmented Lagrangian gradient \eqref{eqn:AugLagGrad} can be evaluated in parallel. On each sub-interval $[t^{\{ \ell \}},t^{\{ \ell+1 \}}]$ in \eqref{eqn:interval-partition} an adjoint solution $ \mathbf{M}^{\rm T}_{\rm k,\rm k+1}\, \lambda_{\rm k+1}$ is computed starting from the terminal value $\lambda_{\rm k+1}$  \eqref{eqn:solution-partition}. The adjoint model runs on each sub-interval can be carried out concurrently. The computational steps are detailed in Algorithm \ref{alg:grad}.
 \begin{algorithm}
\caption{Parallel\_Gradient}\label{alg:grad}
\begin{algorithmic}[1]
\Procedure{Parallel\_Gradient}{}
\State \textbf{Input:} $\lbrack \x_0, \dots, \x_{\rm N};\, \lambda_1, \dots, \lambda_{\rm N} \rbrack$
\State \textbf{Output:} $\lbrack \nabla_{\x_0}\mathcal{L}, \dots, \nabla_{\x_{\rm N}}\mathcal{L}\rbrack$
\ParFor{$0\leq {\rm k} \leq N-1$}
\State $\Delta\x_{\rm k+1} \gets \,\x_{\rm k+1} - \mathcal{M}_{\rm k,k+1}\left(\x_{\rm k}\right)$   \Comment{\parbox[t]{.3\linewidth}{\small{Solution mismatch at sub-interval boundaries}}}
\State $\Delta \y_{\rm k+1} \gets \mathcal{H} \left(\x_{\rm k+1}\right) - \y_{\rm k+1}$
\Comment{\parbox[t]{.3\linewidth}{\small{Innovation vector}}}
\EndParFor 
\ParFor{$0\leq {\rm k} \leq N-1$}
\State $\mathbf{b}_{\rm k+1} \gets  \mu\, \mathbf{P}^{-1}_{\rm k+1}\, \Delta\x_{\rm k+1} - \lambda_{\rm k+1} $   \Comment{\parbox[t]{.3\linewidth}{\small{Translate and scale solution mismatch}}}
\State $\mathbf{d}_{\rm k+1} \gets \mathbf{H}_{\rm N}^{\rm T}\, \R_{\rm N}^{-1}\, \Delta \y_{\rm k+1} $
\Comment{\parbox[t]{.3\linewidth}{\small{Scale innovation vector}}}
\EndParFor 
\ParFor{
$0 \leq \rm k \leq N-1$}
\Comment{\parbox[t]{.3\linewidth}{\small{Perform the adjoint integrations in parallel}}}
\State  $\A_{\rm k} \gets \mathbf{M}^{\rm T}_{\rm k,\rm k+1} \, \mathbf{b}_{\rm k+1}$
\EndParFor
\ParFor{$1\leq \rm k \leq N-1$}   \Comment{\parbox[t]{.3\linewidth}{\small{Compute gradient using \eqref{eqn:AugLagGrad}}}}
\State $\nabla_{\x_{\rm k}}\mathcal{L}\gets  \mathbf{b}_{\rm k} +  \mathbf{d}_{\rm k}-\A_{\rm k}  $
\EndParFor
\State $\nabla_{\x_{\rm 0}}\mathcal{L} \gets \B_0^{-1} \left(\x_0 - \x_0^{\rm b}\right) - \A_0  $ 
\State $\nabla_{\x_{\rm N}} \mathcal{L} \gets \mathbf{d}_{\rm N} + \mathbf{b}_{\rm N}$ 
\EndProcedure
\end{algorithmic}
\end{algorithm}

\subsection{Initial solution guess}

The optimization needs to start with some initial guess for $\x^{\{0\}}$ and $\bm{\lambda}^{\{0\}}$. The initial guess for the state is obtained by performing a serial forward integration using the background value for initial conditions, $\x^{\{0\}}_{\rm k} = \x^{\rm b}_{\rm k}$. 
The initial value for the adjoint variable could be obtained by running once the adjoint model once serially along the background trajectory. In our experiments we choose the simpler, and less expensive, initialization $\lambda^{\{0\}}_{\rm k}=0$.

%
\subsection{Updating Lagrange multipliers}

In order to accelerate the convergence of the optimization process we replace the standard updates \eqref{eqn:update} with the strategy proposed in \cite{he2010acceleration}. The new update process uses information from the previous two iterations instead of just one iteration. The process takes the following steps:
\begin{enumerate}
 \item Choose $\bm \lambda^0$ and set $t^{\{1\}} = 1$.
 \item Let $\x^{\rm \{\ell\}}$ be the solution obtained by solving the optimization problem \eqref{eqn:ipAug} for particular values of $\mu^{\rm \{\ell\}}$ and $\bm\lambda^{\{\ell\}}$. Apply the classical update \eqref{eqn:update} to obtain $\mu^{\{\ell+1\}}$ and $\widetilde\lambda^{\{\ell+1\}}$.
 \item The updated Lagrange multiplier is obtained as follows: 
 \begin{eqnarray*}
\displaystyle t^{\{\ell+1\}} &=&\frac{1}{2}\, \left(1 + \sqrt{1+4 \, (t^{\{\ell\}})^2}\right), \\
\displaystyle \bm\lambda^{\{\ell+1\}} &=& \widetilde {\bm \lambda}^{\{\ell+1\}} + \left(\frac{t^{\{ \ell \}} -1}{t^{\{ \ell+1 \}}}\right) \left( \widetilde{\bm \lambda}^{\{\ell +1\}} - \widetilde {\bm \lambda}^{\rm \{\ell\}} \right) + \left(\frac{t^{\{ \ell \}}}{t^{\{ \ell+1 \}}}\right) \left(\widetilde{\bm \lambda}^{\{\ell +1\}} - \bm \lambda^{\rm \{\ell\}}\right)\,.
 \end{eqnarray*}
\end{enumerate}
It is important to note that above procedure requires the values of $\widetilde{\bm\lambda}$ from two successive outer iterations, namely $\ell$ and $\ell+1$.

\section{Numerical experiments} \label{sec:exp}

We study the performance of the parallel implementation of augmented Lagrangian 4D-Var algorithm using the Lorenz-96 model with 40 variables \cite{lorenz1996} and a shallow water on the sphere model with $\sim$8,000 variables \cite{Amik:2007}. 

\subsection{Lorenz-96 model}
The Lorenz-96 model \cite{lorenz1996} is given by 
      \begin{equation}            
      \label{eqn:Lorenz96}
               \frac{\mathrm{d}\x_{\rm k}}{\mathrm{d}t} = \x_{\rm k-1} \left( \x_{\rm k+1} - \x_{\rm k-2} \right) - \x_{\rm k} + F \,,
               \quad k=1,\dots,40,
      \end{equation}
with periodic boundary conditions and the forcing term $F=8$ \cite{lorenz1996}. We use synthetic observations  generated by perturbing the reference trajectory with normal noise with mean zero and standard deviation of $5\%$ of average magnitude of the reference solution. The background uncertainty is set to $8\%$ of average magnitude of the reference solution. The background and observation error covariance matrices are assumed to be diagonal. A vector of equidistant components ranging from $-2$ to $2$ was integrated forward in time for 200 time steps and the final state is taken as a reference initial condition for the experiments. 
      
\subsection{Shallow water model on the sphere}
The shallow water equations have been used extensively as a simple model of the atmosphere since they contain the essential wave propagation mechanisms found in general circulation models  \cite{Amik:2007}. The shallow water equations in spherical coordinates are:
\begin{subequations}
\label{eqn:swe}
\begin{eqnarray}
 \frac{\partial u}{\partial t} + \frac{1}{a\cos \theta} \left( u \frac{\partial u}{\partial \lambda} + v \cos \theta \frac{\partial u}{\partial \theta} \right) - \left(f + \frac{u \tan \theta}{a} \right) v + \frac{g}{a \cos \theta} \frac{\partial h} {\partial \lambda} = 0, \\
 \frac{\partial v}{\partial t} + \frac{1}{a\cos \theta} \left( u \frac{\partial v}{\partial \lambda} + v \cos \theta \frac{\partial v}{\partial \theta} \right) + \left(f + \frac{u \tan \theta}{a} \right) u + \frac{g}{a} \frac{\partial h} {\partial \theta} = 0, \\
 \frac{\partial h}{\partial t} + \frac{1}{a \cos \theta} \left(\frac{\partial\left(hu\right)}{\partial \lambda} + \frac{\partial{\left(hv \cos \theta \right)}}{\partial \theta} \right) = 0.
\end{eqnarray}
\end{subequations}
Here $f$ is the Coriolis parameter given by $f = 2 \Omega \sin \theta$, where  $\Omega$ is the angular speed of the rotation of the Earth, $h$ is the height of the homogeneous atmosphere, $u$ and $v$ are the zonal and meridional wind components, respectively, $\theta$ and $\lambda$ are the latitudinal and longitudinal directions, respectively, $a$ is the radius of the earth and $g$ is the gravitational constant. The space discretization is performed using the unstaggered Turkel-Zwas scheme \cite{Navon19911311, navon1987application}. The discretized spherical grid has nlon=36 nodes in longitudinal direction and nlat=72 nodes in the latitudinal direction. The semi-discretization in space leads to a discrete model of the form \eqref{eqn:model}.
 In \eqref{eqn:model} the zonal wind, meridional wind and the height variables are combined into the vector $\x \in \mathbb{R}^n$ with $n=3\times{\rm nlat}\times{\rm nlon}$. We perform the time integration using an adaptive time-stepping algorithm. For a tolerance of $\displaystyle 10^{-8}$ the average time step size of the time-integrator is 180 seconds. A reference initial condition is used to generate a reference trajectory.

 Synthetic observation errors at various times $t_{\rm k}$ are normally distributed  with mean zero and a diagonal observation error covariance matrix with entries equal to $(\mathbf{R}_{\rm k})_{\rm i,i}=1$ for $u$ and $v$ components and
 $(\mathbf{R}_{\rm k})_{\rm i,i}=10^6$ for $h$ components. The $\mathbf{R}_{\rm k}$ values correspond to a standard deviation of $5\%$ for $u$ and $v$ components, and $2\%$ for $h$ component. We construct a flow dependent background error covariance matrix as described in \cite{Attia:2014, Attia:2015}. The standard deviation of the background errors for the height component is 2\% of the average magnitude of the reference height component in the reference initial condition. The standard deviation of the background errors for the wind components is 15\% of the average magnitude of the reference wind component in the reference initial condition.
      
%
\subsection{Experimental setup}\label{subsec:Tests}
All numerical experiments are carried out in {\sc Matlab}. The parallel implementations are developed using {\sc Matlab}'s parallel toolbox. 
We compare the performance of the proposed parallel implementation with that of the standard 4D-Var. The accuracy of numerical solutions is measured by the root mean square error (RMSE) with respect to a reference solution. The RMSE is given by
            \begin{equation}
            \label{eqn:RMSE}
                  \textbf{RMSE} = \sqrt{\frac{1}{N}\, \sum_{k=1}^N \nvar^{-1} \, \left\Vert \x^{\rm a}_{\rm k} -\x^{\rm reference}_{\rm k} \right\Vert^2} \,,
            \end{equation}
where the reference solution $\x^{\rm reference}$ and the analysis $\x^{\rm a}$ are propagated forward in time using the full model, and the difference is measured at all times throughout the assimilation window.

Computing the cost functions and gradients is the most important aspect of the 4D-Var algorithm and hence it is necessary that their computations are scalable. To evaluate the gradient and cost function we carry out numerical integrations of the forward and adjoint models using MATLODE \cite{adaug:2015}. MATLODE is a {\sc Matlab} version of FATODE, which was developed in {\sc Fortran} \cite{zhang2014fatode}. The optimization is carried out using the L-BFGS-B solver \cite{liu1989limited} implemented in the Poblano optimization toolbox developed at Sandia National Laboratory \cite{SAND2010-1422}. 

            \subsection{Results with the Lorenz-96 model}
\begin{figure}
  \centering
  {\includegraphics[width=0.8\linewidth]{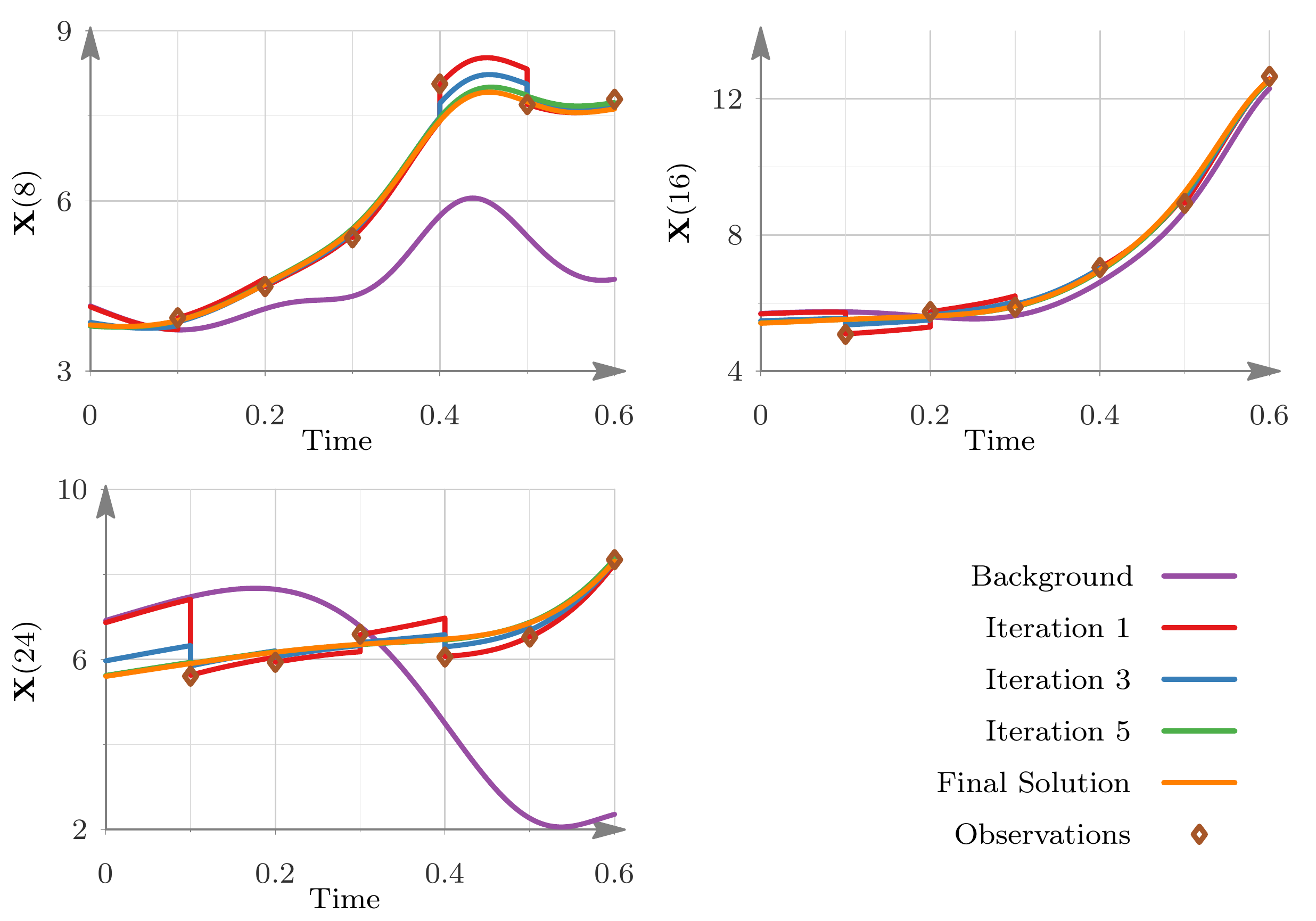}}
  \caption{Lorenz-96 model \eqref{eqn:Lorenz96} results. Convergence of the parallel augmented Lagrangian 4D-Var solution over the first several iterations. The final solution is computed with the traditional 4D-Var. Three different variables are shown as an illustration.}
  \label{fig:Iterations_Lorenz}
\end{figure}

\begin{figure}
  \centering
  {\includegraphics[width=0.8\linewidth]{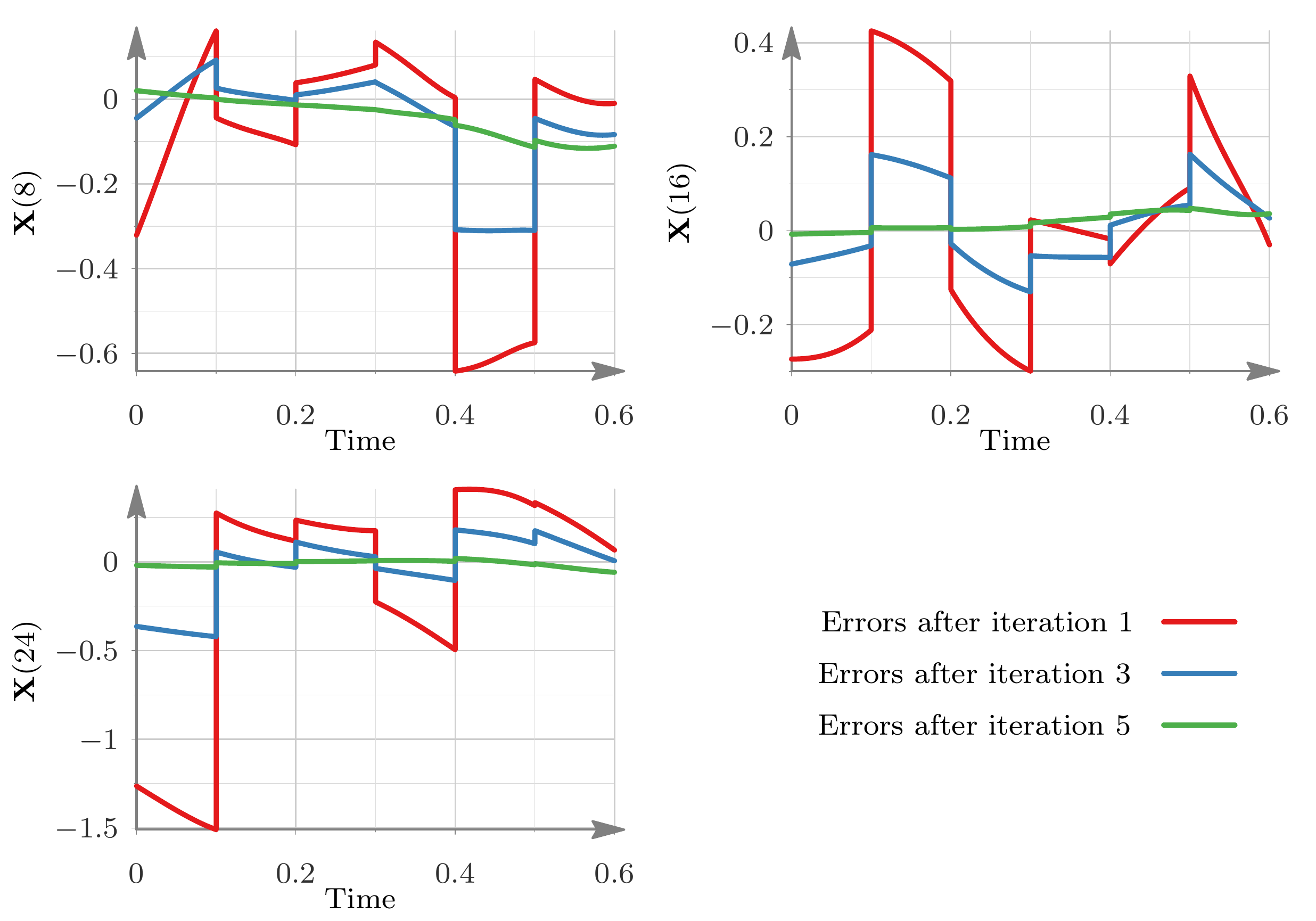}}
  \caption{Lorenz-96 model \eqref{eqn:Lorenz96} results. Difference between the augmented Lagrangian solutions at different iterations and the solution of traditional 4D-Var. Three different variables are shown as an illustration.}
  \label{fig:Errors_Lorenz}
\end{figure}

Figure \ref{fig:Iterations_Lorenz} illustrates the convergence of augmented Lagrangian 4D-Var iterations. The intermediate solutions are discontinuous at sub-interval boundaries. The corresponding errors, with respect to the traditional 4D-Var solution, are shown in Figure \ref{fig:Errors_Lorenz}. They decrease quickly to zero, showing that the augmented Lagrangian 4D-Var solution converges to the solution of traditional 4D-Var.

\begin{figure}
  \centering
  {\includegraphics[width=0.45\linewidth]{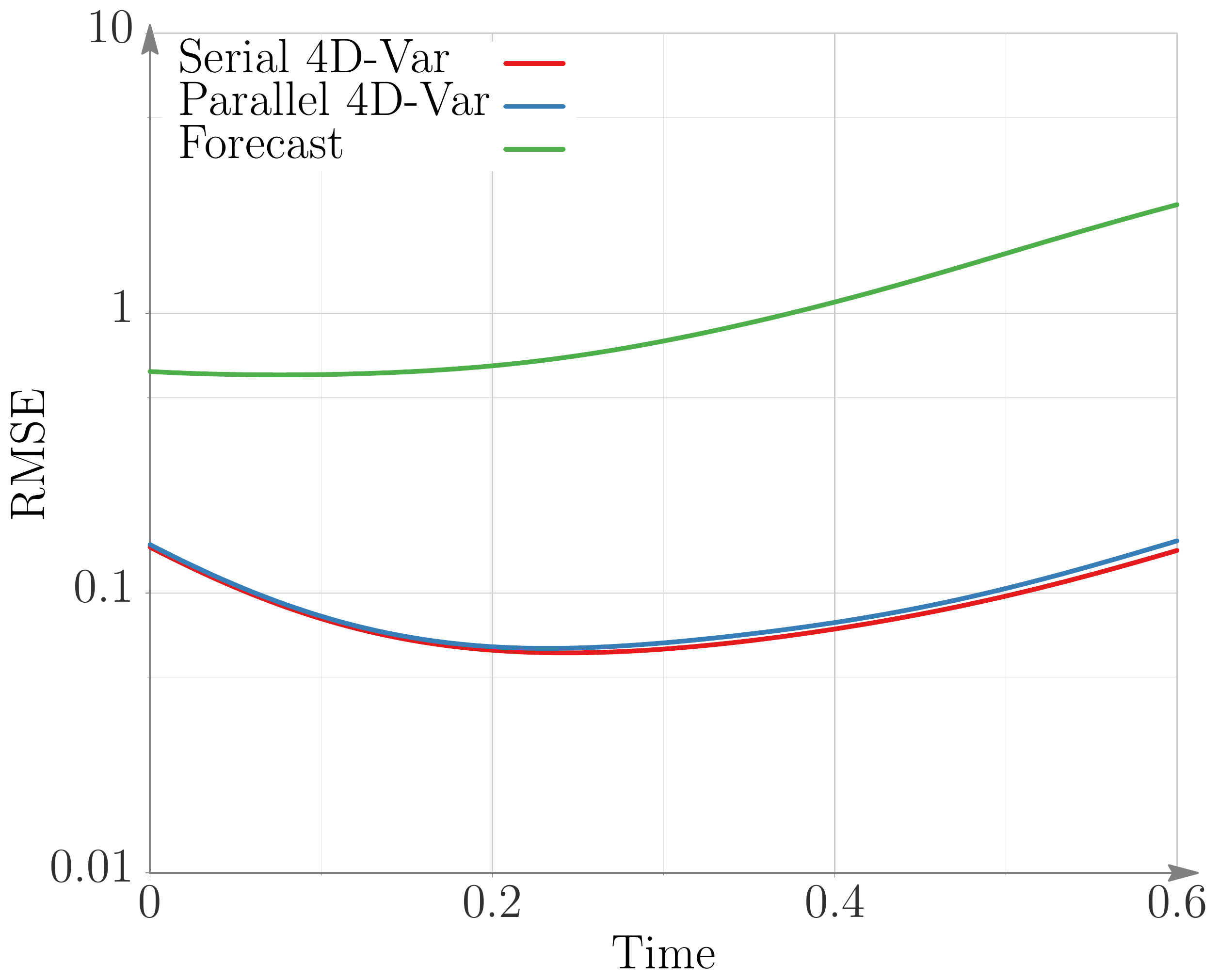}}
  \caption{Lorenz-96 model \eqref{eqn:Lorenz96} results. RMSE errors of traditional and augmented Lagrangian solutions with respect to the reference analysis. The two implementations give nearly identical results. Traditional 4D-Var approach requires 230 iterations, whereas the augmented Lagrangian requires 100 iterations.}
  \label{fig:RMSE_Lorenz}
\end{figure}

\begin{figure}
  \centering
  \includegraphics[width=0.45\linewidth]{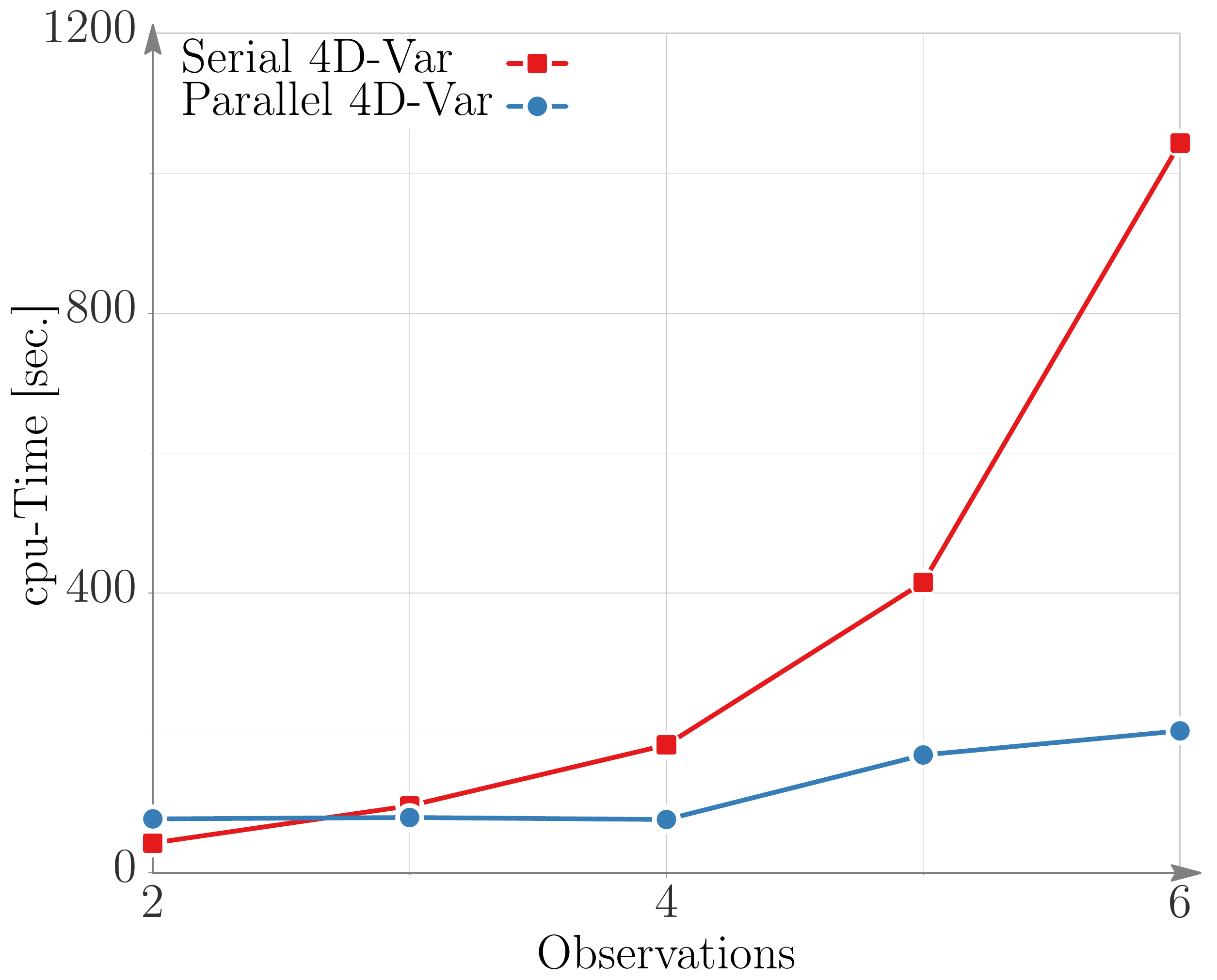}
  \caption{Lorenz-96 model \eqref{eqn:Lorenz96} results. The scaling of overall computational times with increasing number of sub-intervals (which leads to an increase in the length of the assimilation window).}
  \label{fig:Cputime_Lorenz}
\end{figure}

\begin{figure}
  \centering
  \subfigure[Cost of one cost function evaluation for increasing number of sub-intervals\label{fig:CfEvaluation}]{\includegraphics[width=0.45\linewidth]{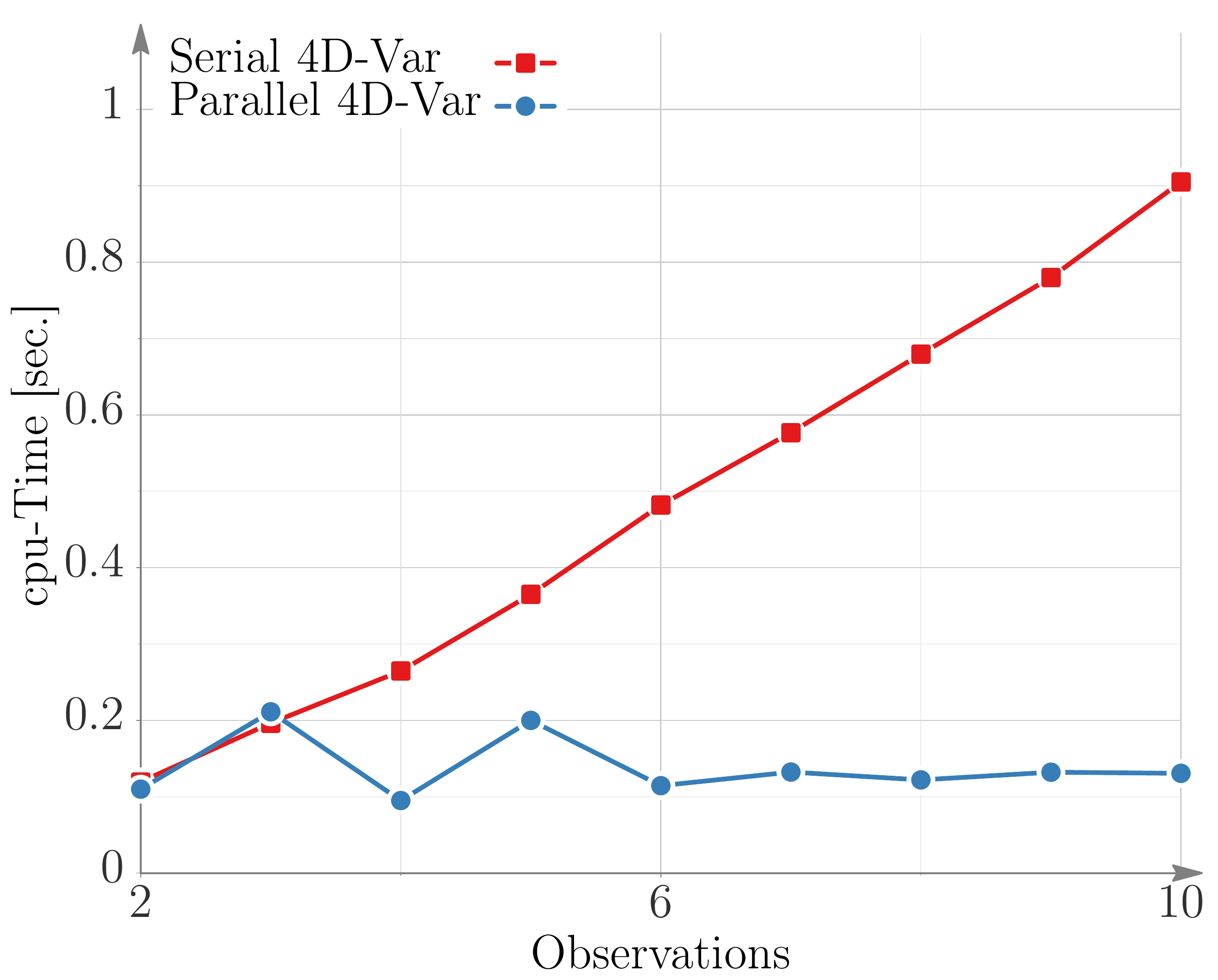}}
  \subfigure[Cost of one cost gradient evaluation for increasing number of sub-intervals\label{fig:GfEvaluation}]{\includegraphics[width=0.45\linewidth]{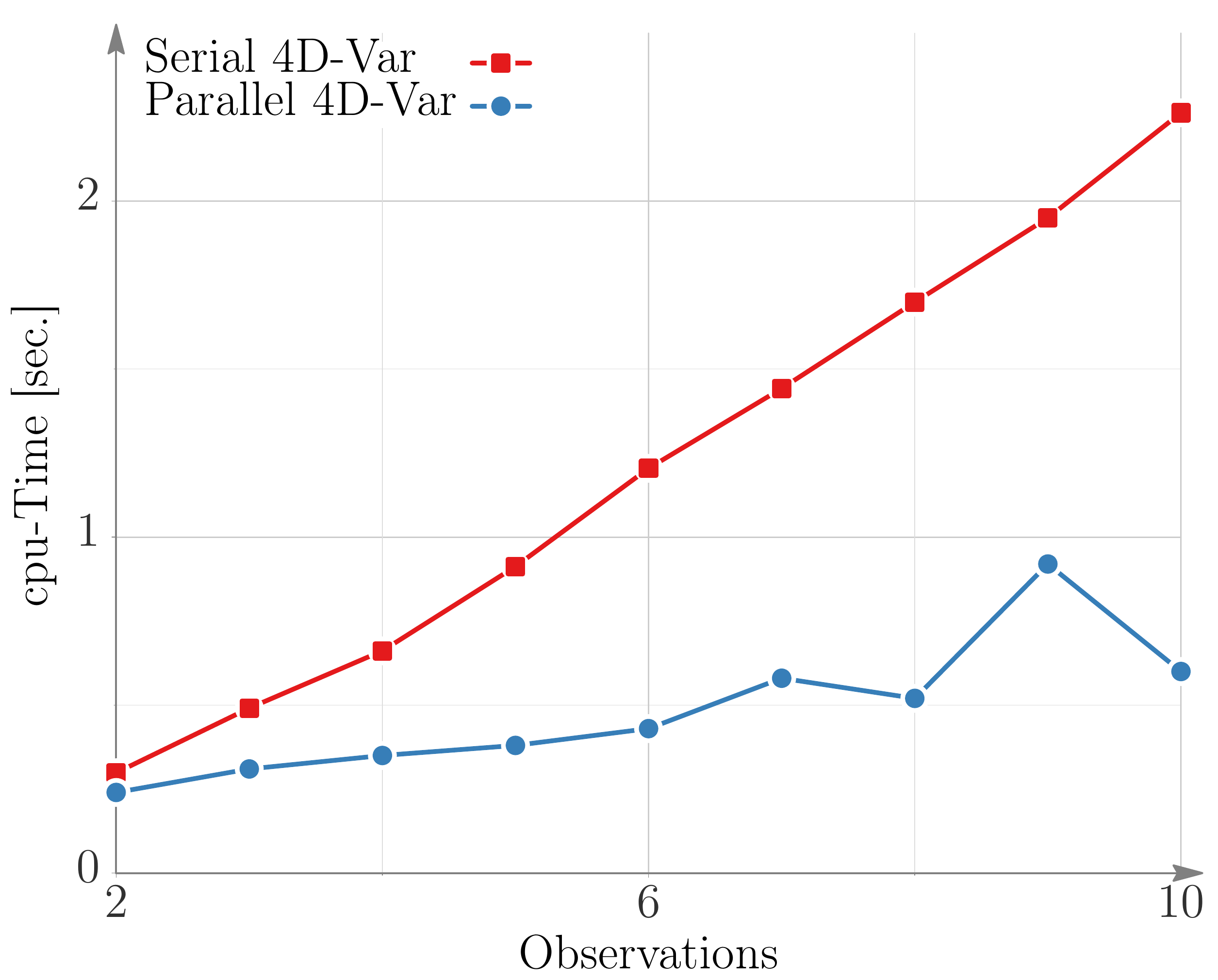}}
  \caption{Lorenz-96 model \eqref{eqn:Lorenz96} results. Weak scalability results for the cost function and gradient evaluations.}
   \label{fig:CfGf_Lorenz}
\end{figure}

The errors of the augmented Lagrangian and traditional 4D-Var solutions with respect to the reference solution are shown in Figure \ref{fig:RMSE_Lorenz}. The reference solution is obtained by propagating the reference initial condition using the forward model in \eqref{eqn:Lorenz96}. The sequential 4D-Var requires $230$ gradient and $574$ cost function evaluations, where as the parallel 4D-Var requires a total of $100$ gradient and $650$ cost function evaluations for $6$ observations. The Weak scalability results are presented in Figure \ref{fig:Cputime_Lorenz}. As the length of the assimilation window increases, the number of sub-intervals increases, and so does the number of observations (one per sub-interval). The number of cores on which the parallel algorithm is run increases such as to remain equal to the number of sub-intervals. The parallel algorithm is scalable in weak sense: the total computational time increases very slowly with an increased problem size. The most time consuming calculations are those of the cost function and gradient evaluations, which require running the forward and the adjoint models, respectively. The results shown in Figure \ref{fig:CfGf_Lorenz} confirm the good weak scalability of the cost function and gradient computations. 

%

 \subsection{Results with the the shallow water model}
\begin{figure}
  \centering
\subfigure[Cost of one cost function evaluation for increasing number of sub-intervals\label{fig:CfEvaluation_SWE}]{\includegraphics[width=0.45\linewidth]{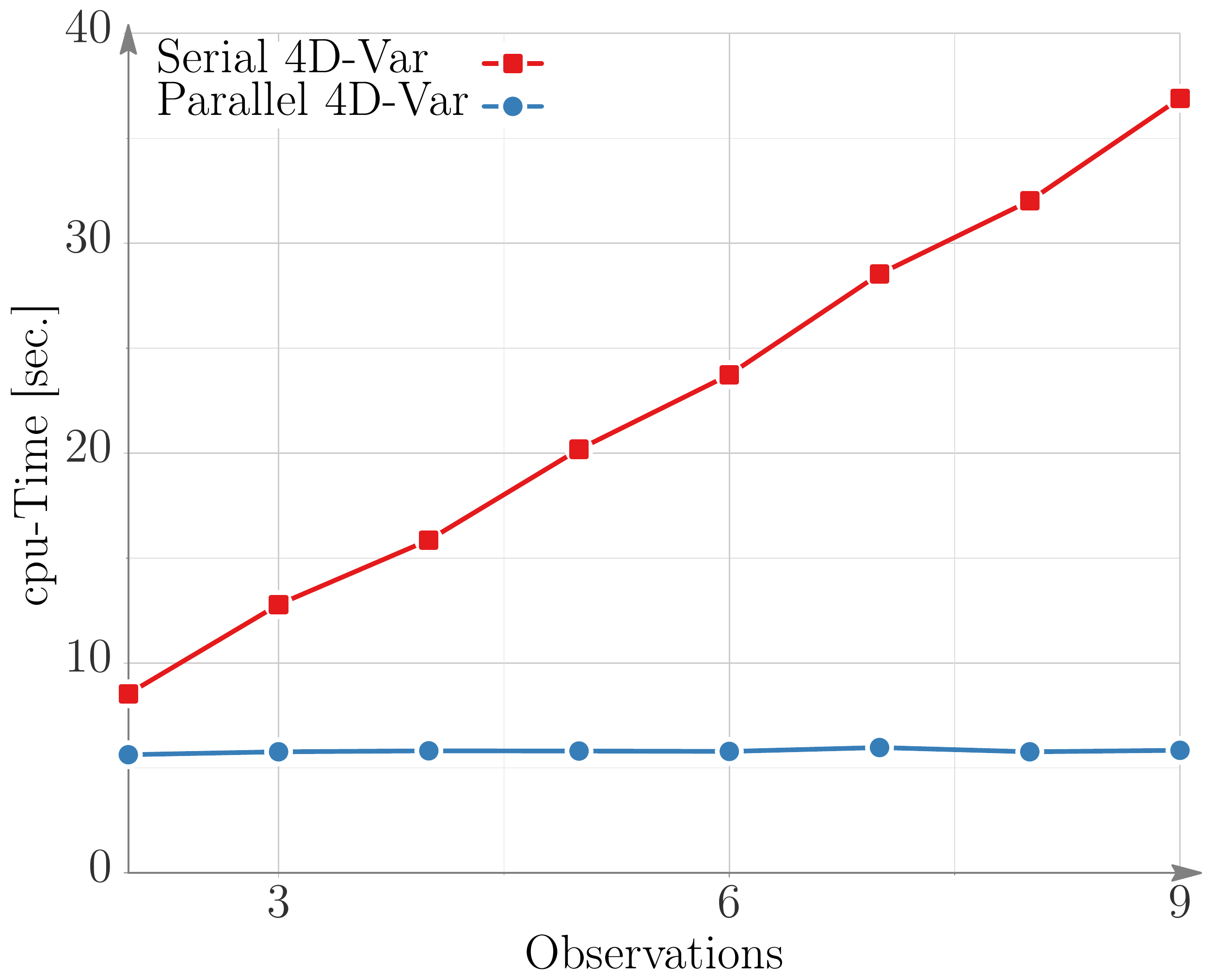}}
\subfigure[Cost of one gradient evaluation for increasing number of sub-intervals\label{fig:GfEvaluation_SWE}]{\includegraphics[width=0.45\linewidth]{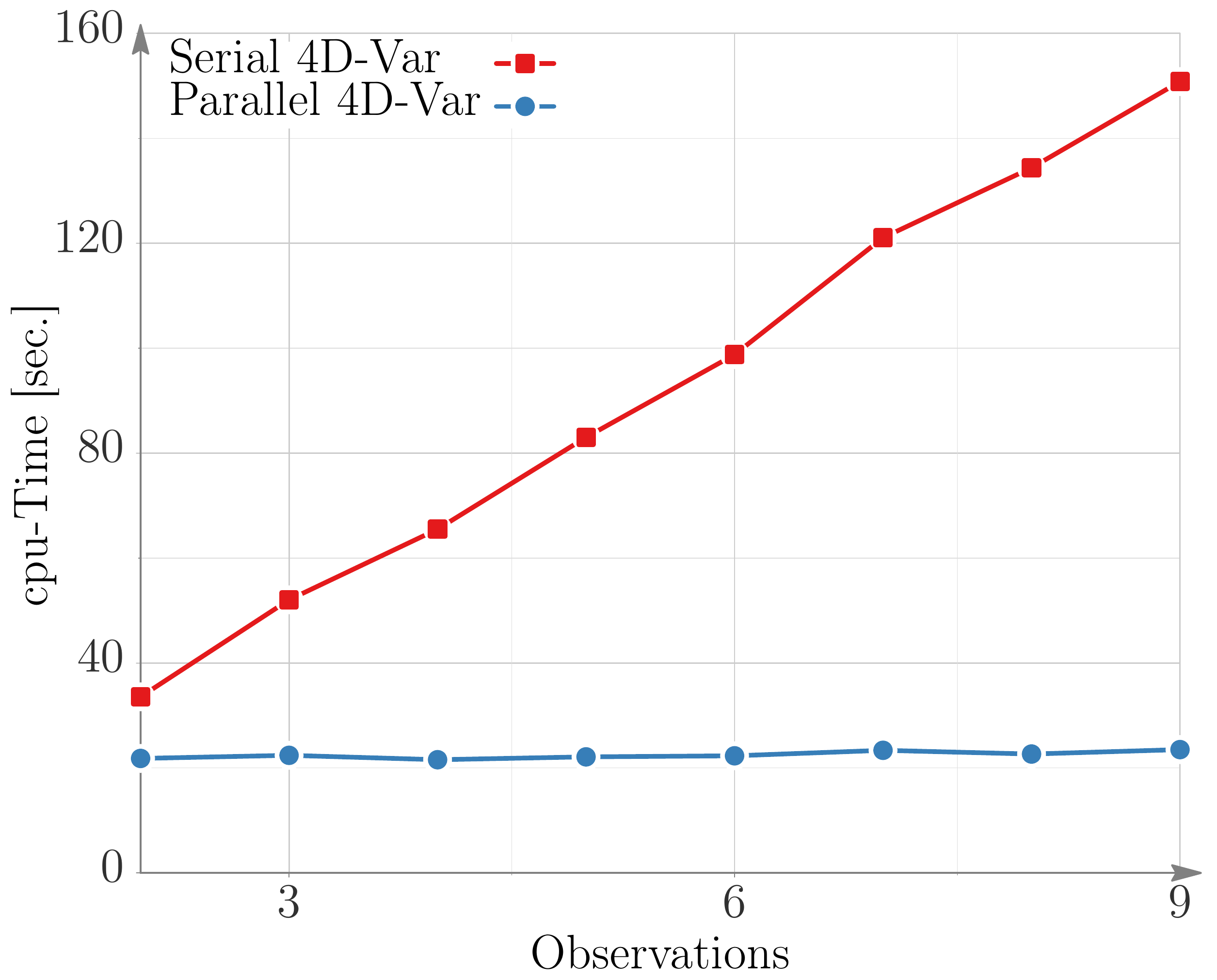}}
  \caption{Shallow water equations \eqref{eqn:swe} results. Weak scalability results for the cost function and gradient evaluations.}
\label{fig:CfGf_SWE}
\end{figure}

\begin{figure}
  \centering
  \subfigure[Hybrid, parallel and serial (traditional) 4D-Var\label{fig:WorkPerformanceHybrid}]{\includegraphics[width=0.45\linewidth]{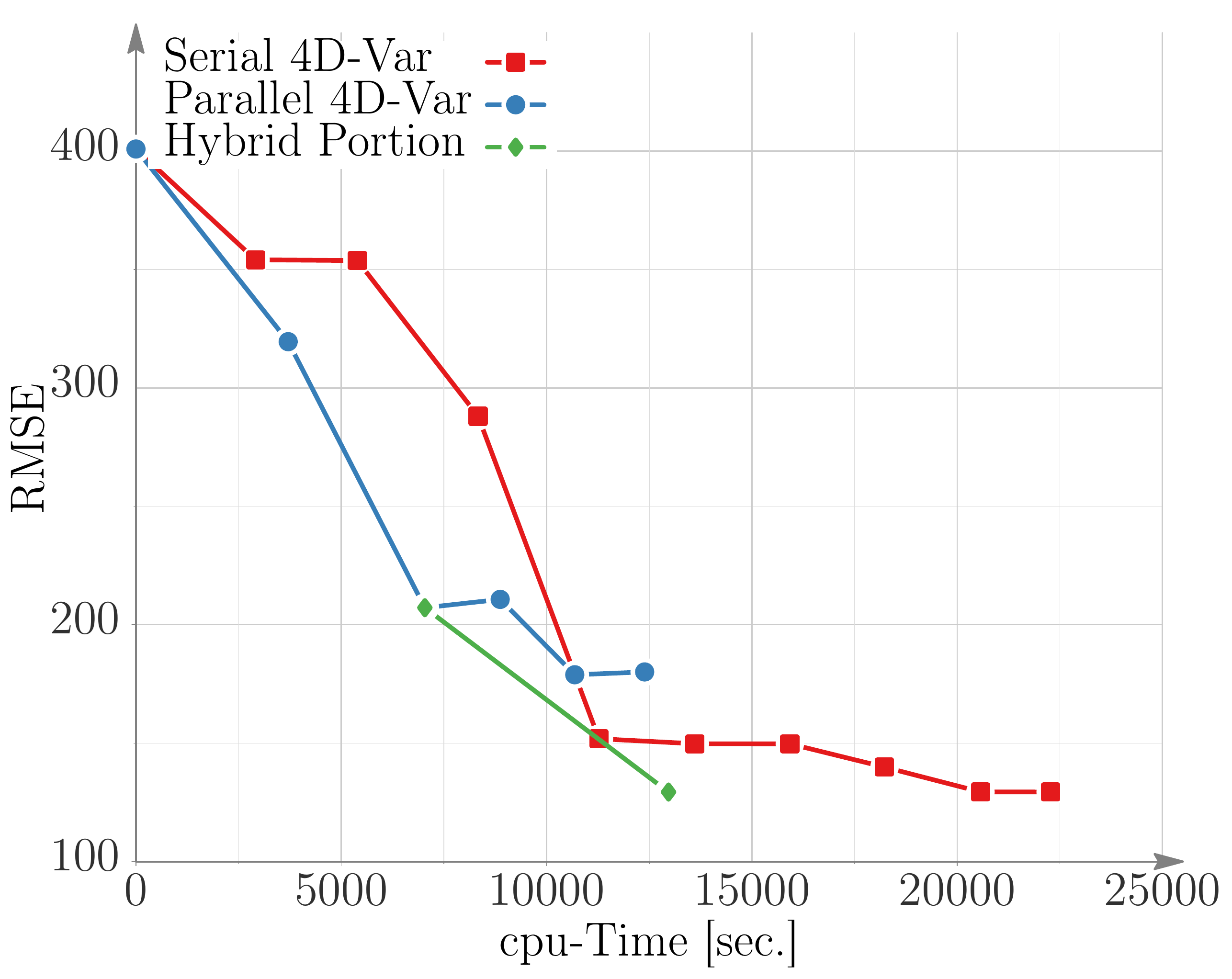}}
  \caption{Shallow water equations \eqref{eqn:swe} results. Work-performance diagrams comparing parallel and hybrid methods with serial 4D-Var for 9 sub-intervals.}
  \label{fig:WorkPerformance_SWE}
\end{figure}

Figure \ref{fig:CfGf_SWE} shows the weak scalability of cost function and gradient evaluations with the shallow water model. It can be seen that in both cases the parallel computational time is nearly constant with increasing problem size (number of sub-intervals) and a proportional increase in the number of cores (the number of cores is equal to the number of sub-intervals). 
Figure \ref{fig:WorkPerformance_SWE} presents the work-precision diagrams, i.e., the evolution of solution accuracy (RMSE) with increasing number of iterations (increasing CPU time). The initial iterates of the parallel 4D-Var solution proceed rapidly, but afterwards the convergence slows down. At this stage the penalty parameter $\mu$ fairly large and the optimization problem becomes more difficult for the LBFGS algorithm. The performance can be improved by replacing LBFGS with algorithms specially tailored to solve optimization problems in the augmented Lagrangian framework \cite{birgin2014practical}. An alternative strategy is to use a hybrid method: employ parallel 4D-Var for several iterations in the beginning, then continue with traditional serial 4D-Var. Here we perform two outer iterations with small values of $\mu$ and then use this solution as an initial guess for the serial 4D-Var. This strategy improves the performance considerably as seen in Figure \ref{fig:WorkPerformanceHybrid}. Table \ref{tab:timings} provides the computational times for parallel and serial 4D-Var algorithms. The serial 4D-Var for $9$ observations requires $95$ gradient and $200$ cost function evaluations, whereas the parallel 4D-Var algorithm requires $350$ gradient and $720$ cost function evaluations. In the hybrid methodology, we perform $200$ iterations of parallel 4D-Var and $30$ iterations of serial 4D-Var for $9$ observations. The final RMSE over the assimilation window for both the serial 4D-Var and hybrid methods is $\sim$ 125. We notice a steady increase in speedup as the problem size (number of sub-intervals) is increased. It is possible to further improve the performance of the parallel algorithm by using second derivative information in the form of Hessian-vector products and employ Newton-type methods \cite{cioaca2012second}. 
\begin{table}[ht]
 \centering
 \begin{tabular}{|c|c|c|c|}
  \hline
  No. of sub-intervals & Serial time [sec.] & Parallel/Hybrid time [sec.] & Speedup  \\ \hline
  5 & 8,515 & 10,834 & 0.786\\ \hline
  7 & 15,745 & 12,782 & 1.232\\ \hline
  9 & 22,277 & 12,971 & 1.756\\ \hline
 \end{tabular}
 \caption{Shallow water equations \eqref{eqn:swe} results. Performance comparison of augmented Lagrangian/parallel and traditional/serial 4D-Var and the corresponding speedups.}
 \label{tab:timings}
\end{table}
%
\section{Conclusions and future work}\label{sec:conc}
This work presents an augmented Lagrangian framework to perform strong-constraint 4D-Var data assimilation in parallel. The assimilation window is split in sub-intervals; cost function and gradient evaluations, which are the main components of the algorithm, are performed by running the forward and the adjoint model in parallel across different sub-intervals. To the best of our knowledge this is the first work that uses an augmented Lagrangian approach to data assimilation, and the first one to propose a time-parallel implementation of strong-constraint 4D-var.

Future work will focus on tuning the optimization procedure to improve performance on large scale problems, e.g.,  data assimilation with the weather research and forecast model \cite{skamarock2005description}. The size of the control variable in the augmented Lagrangian framework increases with the number of sub-intervals and can become a bottleneck for the optimization. One possible strategy to overcome this is to perform optimization on a coarse grid and use the projected solution as an initial guess for the fine grid. A natural extension to our methodology is to use the augmented Lagrangian framework to expose space-time parallelism in the data assimilation problem such as to create more parallel tasks and improve the overall scalability. Space parallelism in a penalty formulation has been recently discussed in \cite{Damore:2014}. We will consider the use of optimization algorithms that are specifically tuned to work well in the augmented Lagrangian framework \cite{birgin2014practical}. Next, Hessian information can be used to accelerate the convergence significantly when the iterates are close to minima. In order to implement this it is useful to explore the construction of  second  order adjoint models that compute the Hessian-vector products in parallel.

\section*{Acknowledgements}
This work was supported in part by the awards AFOSR FA9550--12--1--0293--DEF, AFOSR 12-2640-06, NSF DMS--1419003, NSF CCF--1218454, and by the Computational Science Laboratory at Virginia Tech.

\section*{References}


\end{document}